\documentclass[12pt]{amsart}
\usepackage{euscript}
\usepackage{amssymb,amscd}
\theoremstyle{plain}
\newtheorem{theorem}{Theorem}[section]

\newtheorem{lemma}[theorem]{Lemma}
\newtheorem{proposition}[theorem]{Proposition}
\newtheorem{corollary}[theorem]{Corollary}

\newtheorem*{thma}{Theorem A}
\newtheorem*{thmb}{Theorem B}
\newtheorem*{thmc}{Theorem C}
\newtheorem*{thmd}{Theorem D}
\newtheorem*{thme}{Theorem E}
\newtheorem*{thmf}{Theorem F}
\newtheorem*{thmg}{Theorem G}
\newtheorem*{thmh}{Theorem H}
\newtheorem*{quest}{Question}

\newtheorem{definition}[theorem]{Definition}
\newtheorem{example}[theorem]{Example}
\newtheorem{remark}[theorem]{Remark}

\numberwithin{equation}{section}

\newcommand{\R}{\mathbb{R}}
\newcommand{\RR}{\overline{\mathbb{R}}_+}
\newcommand{\Z}{\mathbb{Z}}
\newcommand{\NN}{\mathbb{N}_0}
\newcommand{\N}{\mathbb{N}}
\newcommand{\eps}{\varepsilon}
\renewcommand{\phi}{\varphi}
\renewcommand{\kappa}{\varkappa}

\begin{document}
\title[Affine actions of a free semigroup on the real line]
{Affine actions of a free semigroup\\ on the real line}

\author{Vitaly Bergelson}
\address{Department of Mathematics, Ohio State University, Columbus,
OH 43210}
\email{vitaly@math.ohio-state.edu}

\author{Micha{\l} Misiurewicz}
\address{Department of Mathematical Sciences,
Indiana University Purdue University Indianapolis, Indianapolis, IN
46202-3216}
\email{mmisiure@math.iupui.edu}

\author{Samuel Senti}
\address{Instituto de Matematica Pura e Aplicada,
Rio de Janeiro 22460-320, Brazil}
\email{senti@impa.br}

\date{December 17, 2005}

\thanks{The first author was partially supported by NSF grant DMS
  0345350, the second author by NSF grant DMS 0456526, and the third
  author by a grant from the Swiss National Science Foundation. The
  authors are grateful to Mariusz Urba\'nski and Boris Solomyak for
  comments on the first version of the paper.}

\begin{abstract}
We consider actions of the free semigroup with two generators on the
real line, where the generators act as affine maps, one contracting and one
expanding, with distinct fixed points. Then every orbit is dense in a
half-line, which leads to the question whether it is, in some sense, uniformly
distributed. We present answers to this question for various
interpretations of the phrase ``uniformly distributed''.
\end{abstract}

\maketitle


\section{Introduction}\label{intro}


Denote the set of all one-sided 0-1 sequences by
$\Sigma=\{0,1\}^{\NN}$, where $\NN=\{0,1,2,\dots\}$, and let
$\Sigma_n$ be the set of all $0$-$1$ sequences of length $n$
($\Sigma_0$ consists of one empty sequence). Set
$$G=\bigcup_{n\in\NN}\Sigma_n.$$
When the semigroup multiplication is the concatenation, then $G$ is a
free semigroup with two generators.

Suppose that for some space $Y$ two maps, $T_0,T_1:Y\to Y$, are given.
For $\omega=(\omega_0,\omega_1,\dots,\omega_{n-1})\in \Sigma_n$ write
$$T_{\omega}=T_{\omega_{n-1}}\circ \dots\circ T_{\omega_1}\circ
T_{\omega_0}.$$
Then $G$ acts on $Y$ by $\omega\mapsto T_\omega$.

Our starting point is the following result, equivalent to
Corollary~3.2 of \cite{Mis-Rod1}.

\begin{thma}
Let $T_0(x)=x/2$ and $T_1(x)=(3x+1)/2$, $x\in\R_+=[0,\infty)$. Then
for any $x\in\R_+$ the orbit $\{T_\omega(x)\}_{\omega\in G}$ is dense
in $\R_+$.
\end{thma}

We now consider a somewhat more general setup, where we have one expanding and
one contracting affine transformation of $\R_+$ with distinct fixed
points. Any pair of such transformations may be
brought into the form $T_0(x)=ax$ and $T_1(x)=bx+1$, where $0<a<1<b$,
by conjugating via an affine map. Therefore we will assume from now on
that $T_0$ and $T_1$ have such form.

Both Theorem~3.1 and Theorem~3.4 of \cite{Mis-Rod1} imply Theorem~A.
However, the proof of Theorem~3.1 of \cite{Mis-Rod1} cannot be
generalized for all $T_0$ and $T_1$ as defined in the preceding
paragraph, since it requires rational independence of $\log a$ and
$\log b$. Fortunately, one can
repeat almost verbatim the proofs of Lemma~3.3 and Theorem~3.4 of
\cite{Mis-Rod1}, which gives the following result.

\begin{thmb}
Let $T_0(x)=ax$ and $T_1(x)=bx+1$, $x\in\R_+=[0,\infty)$, where
$0<a<1<b$. Then for any $x\in\R_+$ the orbit
$\{T_\omega(x)\}_{\omega\in G}$ is dense in $\R_+$.
\end{thmb}

Theorem~B naturally leads to the following question.

\begin{quest}
Is it true that under the assumptions of Theorem~B every orbit is
uniformly distributed in $\R_+$?
\end{quest}

One of the problems with this Question is, of course, the vagueness of
the phrase ``uniformly distributed''. Indeed, while in the classical
theory of the uniform distribution of sequences modulo one
(see, e.g., \cite{Kuipers-Nie}) the sequences are conveniently indexed by the
ordered set $\N=\{1,2,3,\dots\}$ and belong to the ``nice'' compact
space ${\mathbb T}=\R/\Z$, none of
these ingredients is present in our situation.

As we shall see, there are various natural approaches to interpreting
and answering our Question, each revealing an interesting facet of the
situation at hand.

One possible approach would be to inquire whether one can order the
elements of our semigroup $G$ in some way, say $(g_i)_{i\in\N}$,
so that for some probability measure\footnote{By a probability
measure we will always understand a Borel probability measure.} $\mu$
on $\RR=[0,\infty]$, one has
\begin{equation}\label{conver}
\lim_{n\to\infty}\frac1n\sum_{i=1}^n f(T_{g_i}(x))=\int_{\RR}f\;d\mu
\end{equation}
for every $x\in\R_+$ and every $f\in C(\RR)$, where $C(\RR)$ is the space of
all continuous real functions on $\RR$. While a priori it is not clear at all
whether one can find such ordering of $G$ and such measure $\mu$, the
following theorem shows that, actually, a much stronger result holds.

\begin{thmc}
For \emph{any probability measure} on $\RR$ there exists a sequence
$(g_i)_{i\in\N}$ of elements of $G$ (which can be chosen in such
a way that it exhausts $G$), such that \eqref{conver} holds for any
$x\in\R_+$ and any $f\in C(\RR)$.
\end{thmc}

We will now describe another approach, which is, in our opinion, much
more natural, since it takes into account the structure of the acting
semigroup.

By extending $T_0$ and $T_1$ to continuous maps of the one-point
compactification $\RR$ of $\R_+$ by the rule
$T_0(\infty)=T_1(\infty)=\infty$ (and keeping the notation $T_0,T_1$
for these extensions), we obtain an action of $G$ on a compact space
$\RR$.
Let $\mathcal{M}$ (respectively $\overline{\mathcal{M}}$) denote the
space of all probability measures on $\R_+$ (respectively $\RR$).

A possible approach to answering our Question is to look for a measure
$\mu\in\overline{\mathcal{M}}$ such that the natural from the point of
view of ergodic theory averages
$$\frac{1}{|\Sigma_n|}\sum_{\omega\in\Sigma_n}f(T_\omega(x))$$
converge to $\displaystyle{\int_{\RR}f\;d\mu}$ for every $f\in C(\RR)$
and every $x\in\R_+$.

A priori, it is not clear whether such measure $\mu$ exists. A
possible candidate for such $\mu$ would be a measure that is
$g$-invariant for any $g\in G$. However,
it is easy to see that there are no such nontrivial measures
in our situation. Indeed, since $T_1(x)>x$ for every $x\in\R_+$,
already the only $T_1$-invariant measure in $\overline{\mathcal{M}}$
is $\delta_\infty$. It therefore comes as a pleasant surprise that
(under an additional assumption on $a$ and $b$) such a measure exists
and can be described quite explicitly.

\begin{thmd}
If, in addition to our standard assumption that $0<a<1<b$, the
parameters $a,b$ satisfy $ab<1$, then there exists $\mu\in\mathcal{M}$
such that for any $f\in C(\RR)$ and any $x\in\R_+$ one has
\begin{equation}\label{con}
\lim_{n\to\infty}\frac{1}{|\Sigma_n|}\sum_{\omega\in\Sigma_n}
f(T_\omega(x))=\int_{\RR} f\;d\mu.
\end{equation}
\end{thmd}

It turns out that when $ab\ge 1$ then the averages in \eqref{con}
converge to $f(\infty)$. However, as we will see, a natural
\emph{weighted} version of Theorem~D holds under much more
general assumptions. In order to give a precise formulation of this
result, we have to introduce additional definitions and notation.

Let $\sigma:\Sigma\to\Sigma$ be the full one-sided shift, i.e.
$\left(\sigma(\bar\omega)\right)_i=\omega_{i+1}$ for any
$\bar\omega=(\omega_i)_{i\in\NN}\in \Sigma$. Set $X=\Sigma\times \RR$
and denote by $\pi_1$ and $\pi_2$ the projections of $X$ onto $\Sigma$
and $\RR$ respectively. Let $\mathcal{N}$ denote the space of all
probability measures on $\Sigma$. Define $\tau:X\to X$ as the skew product
$$\tau(\bar\omega,x)=(\sigma(\bar\omega), T_{\omega_0}x).$$
A similar approach to investigation of actions of free semigroups has
been used for instance by Sumi \cite{Sumi}.

\begin{definition}\label{definv1}
Let $\nu$ be a $\sigma$-invariant probability measure on $\Sigma$. We
will say that a measure $\mu\in\mathcal{M}$ is \emph{invariant for
$(T_0,T_1,\nu)$} if there exists a $\tau$-invariant probability
measure $\kappa$ on $X$ such that the projection $(\pi_1)_*(\kappa)$
to $\mathcal{N}$ is $\nu$ and the projection $(\pi_2)_*(\kappa)$ to
$\overline{\mathcal{M}}$ is $\mu$.
\end{definition}

Note that by choosing in this definition $\mu\in\mathcal{M}$ (rather
than $\mu\in\overline{\mathcal{M}}$) we exclude the case $\mu=\delta_\infty$,
although $\kappa=\nu\times\delta_\infty$ is always $\tau$-invariant
and projects to $\nu$ and $\delta_\infty$.

It makes sense to speak of the Lyapunov exponent of $\nu$, when we
look at the derivative of $\tau$ in the direction in which it exists,
that is, in $\RR$. The derivatives of $T_0$ and $T_1$ are constant and
equal to $a$ and $b$ respectively. We thus define the {\it Lyapunov
exponent} $\lambda(\nu)$ of an ergodic $\sigma$-invariant measure
$\nu\in\mathcal{N}$ as the integral of the function which equals to
$\ln a$ on $C_0$ and $\ln b$ on $C_1$, where
$C_i=\{\omega\in\Sigma:\omega_0=i\}$:
$$\lambda(\nu)=\nu(C_0)\ln a+\nu(C_1)\ln b.$$

\begin{thme}
For any ergodic $\sigma$-invariant probability measure $\nu\in\mathcal{N}$ there
exists a probability measure invariant for $(T_0,T_1,\nu)$ if and only
if $\lambda(\nu)<0$. If such a measure $\mu$ exists, then it is unique
and for every $f\in C(\RR)$ and every $x\in\R_+$ one has
$$\lim_{n\to\infty}\sum_{\omega\in \Sigma_n}\nu(C_\omega)
f(T_\omega(x))=\int_{\RR} f\;d\mu,$$
where for $\omega=(\omega_0,\omega_1,\dots,\omega_{n-1})$,
$C_\omega=\{\bar\omega\in\Sigma:(\bar\omega)_i=\omega_i,
i=0,1,\dots,n-1\}$.
\end{thme}

The issue of what happens (from the point of view of our Question)
when $\lambda(\nu)\ge 0$, remains not fully resolved. One can show
that the weighted averages which appear in Theorem~E converge to
$f(\infty)$. This leads to a suspicion that there is an infinite
measure on $\R_+$ which is ``responsible'' for the statistical
behavior of $G$-orbits. We do not address this issue here.

We also have an interesting addition to Theorem~E, which deals with,
so to say, individual paths of points.

\begin{thmf}
Assume that $\nu\in\mathcal{N}$ is $\sigma$-invariant, ergodic and satisfies
$\lambda(\nu)<0$. Then for $\nu$-almost every
$\bar\omega=(\omega_0,\omega_1,\dots)\in\Sigma$ and for every $f\in
C(\RR)$ and $x\in\R_+$ one has
$$\lim_{n\to\infty}\frac1n\sum_{i=0}^{n-1}f(x_n)=\int_{\RR} f\;d\mu,$$
where $x_n=T_{(\omega_0,\omega_1,\dots,\omega_{n-1})}(x)$ and $\mu$ is
the unique measure invariant for $(T_0,T_1,\nu)$.
\end{thmf}

One more way of addressing our Question is to inquire whether for some
$\sigma$-invariant ergodic measure $\nu\in\mathcal{N}$ the
corresponding $(T_0,T_1,\nu)$-invariant measure $\mu$ is
absolutely continuous with respect to the Lebesgue measure. We prove
the following result.

\begin{thmg}
For every $\gamma>1$ there exists a $\sigma$-invariant ergodic
probability measure $\nu_\gamma\in\mathcal{N}$ with negative Lyapunov
exponent, for which the corresponding $(T_0,T_1,\nu_\gamma)$-invariant
measure $\mu_\gamma$ is absolutely continuous with respect to the
Lebesgue measure and the support of $\mu_\gamma$ is the union of
finitely many intervals whose convex hull is
$[\frac{\gamma-1}{b},\frac{\gamma}{a}]$. In particular,
$\mu_\gamma\neq \mu_\delta$ for $\gamma\neq\delta$, $\gamma,\delta>1$.
\end{thmg}

The case when $\nu$ is $(p,1-p)$ Bernoulli is of a special interest and
needs a separate discussion. If we assume that $a^pb^{1-p}<1$ then
$\lambda(\nu)<0$, and so the $(T_0,T_1,\nu)$-invariant measure $\mu$
exists. If, in addition, $G$ acts effectively (that is,
$\omega\ne\omega'$ implies $T_\omega\ne T_{\omega'}$) and
$a^pb^{1-p}>p^p(1-p)^{1-p}$, then the formal computation of
the Hausdorff dimension of $\mu$, which does not take into account the
overlap of the images of $\R_+$ under $T_0$ and $T_1$, gives
the result larger than 1. In such a case it
is widely conjectured that $\mu$ should be absolutely continuous with
respect to the Lebesgue measure (see, e.g., \cite{PSS}). Unfortunately, the existing
techniques which could be used to prove this, at least for typical
values of $a$ and $b$ (e.g. \cite{SSU}, \cite{Tsujii}, \cite{Rams}),
use strong assumptions which cannot be easily checked in our case. As
a matter of fact, we believe that the above conjecture is false in our situation.

On the other hand, even without assuming that $G$ acts effectively, we
can prove the following interesting property of $\mu$.

\begin{thmh}
Assume that $\nu$ is $(p,1-p)$-Bernoulli with $\lambda(\nu)<0$. Then
the $(T_0,T_1,\nu)$-invariant probability measure $\mu$ on $\R_+$ is the image of
the Lebesgue measure under an increasing map from $[0,1)$ onto $\R_+$,
which is H\"older continuous with the same exponent on each compact
interval.
\end{thmh}

A property of this type is ubiquitous in the Dynamical Systems
Theory. For example, any conjugacy between diffeomorphisms on compact
hyperbolic sets is H\"older continuous (Theorem~19.1.2 of
\cite{Ka-Ha}). Therefore, the measure with maximal
entropy for any Anosov diffeomorphism on a torus is the image of the
Lebesgue measure under a H\"older continuous homeomorphism of the
torus.\footnote{The authors thank Anatole Katok for pointing this out.}

The paper is organized as follows. We deal with Theorem~C in
Section~\ref{semigr}, with existence of $(T_0,T_1,\nu)$-invariant
measures in Section~\ref{exuni}, with Theorems~E and~F in
Section~\ref{erg}, with Theorem~G in Section~\ref{secinv},
and with Theorem~H in Section~\ref{bern}. When stating the results in
those sections, we use a slightly different language, which employs the
notion of the weak-* convergence of measures.


\section{If one disregards generators\dots}\label{semigr}


Let us consider the maps $T_0, T_1$ and the semigroup $G$ defined in
the introduction. It is a well known fact (see for instance
Theorem~2.2 of \cite{Mis-Rod1}) that if $a=1/2$ and $b=3/2$, then the
action of $G$ is effective. In fact, it is easy to see that the set of
$(a,b)$ for which this action is effective contains
the complement of a union of
countably many algebraic curves. Indeed, if the action of $G$ is not
effective, then there exist distinct words $\omega$ and $\omega'$ such
that $T_\omega=T_{\omega'}$. In particular $T_\omega(0)=
T_{\omega'}(0)$ and this is a polynomial equation in $a$, $b$.

\begin{remark}
While rational independence of the logarithms of $a$ and $b$ seems to
play a crucial role for effectiveness of the action of $G$, it is not
sufficient for it. For instance if $a=1/2 $ and $b=4/3$ then
$T_{10001}=T_{01100}$, hence the action is not effective.
\end{remark}

The following theorem (together with Remark~\ref{rem25}) is equivalent
to Theorem~C.

\begin{theorem}\label{approxim}
For any $\mu\in\mathcal{M}$ there exists a infinite sequence
$(g_i)_{i=1}^\infty$ with $g_i\in G$ such that, for any point
$x\in\mathbb{R}_+$ we have
\begin{equation}\label{approxeq}
\frac1n\sum_{i=1}^n\delta_{T_{g_i}x}\stackrel{*}{\longrightarrow}\mu.
\end{equation}
where $\stackrel{*}{\longrightarrow}$ denotes the weak-* convergence
in $\RR$.
\end{theorem}

For the proof of this theorem we will need to adapt Lemma~3.3 and
Theorem~3.4 of \cite{Mis-Rod1} to our more general setup with
$T_0(x)=ax$ and $T_1(x)=bx+1$, with $0<a<1<b$. The proofs are the
same as in \cite{Mis-Rod1}.

\begin{lemma}\label{MR33}
Let $m$ be the number of $1$'s in $\omega=(\omega_0, \ldots,
\omega_{n-1})\in G$ of length $n$. Assume that $x\le M$ and for every $k<n$ we have
$T_{\omega_k}\circ\ldots\circ T_{\omega_0}(x)\le M$. Then
$$\frac{1}{a^{n-m}b^m}\ge (m-1)\frac{x}{bM^2}.$$
\end{lemma}

\begin{theorem}\label{MR34}
Given $x, y>0$ and $\varepsilon>0$, there is $\omega=(\omega_0,
\ldots, \omega_{n-1})\in G$ such that $|T_\omega(x)-y|<\varepsilon$.
Moreover
$$\min(x, y-\varepsilon)\le T_{\omega_k}\circ\ldots\circ
T_{\omega_0}(x)\le K$$
for every $k<n$ where $K$ is a constant depending on $x$ and $y$.
\end{theorem}

\begin{proof}[Proof of Theorem~\ref{approxim}]
For any measure $\mu\in\mathcal{M}$ there exists a sequence
$(a_i)_{i=1}^\infty$ of points $a_i\in\mathbb{R}_+$ such that the
averages of the Dirac delta measures $\delta_{a_i}$ at $a_i$ converge
to $\mu$:
$$\frac1n\sum_{i=1}^n\delta_{a_i}\stackrel{*}{\longrightarrow}\mu.$$
Observe that $\overline{\{T_g(1)\}}_{g\in G}=\R_+$ by
Theorem~\ref{MR34}, and together with Lemma~\ref{MR33}, this implies
that for any $i$ there exists $g_i\in G$ with
\begin{equation}\label{contraction}
|T_{g_i}(1)-a_i|<\frac{1}{i}\qquad\mbox{ and }\qquad
 T_{g_i}'<\frac{1}{i}.
\end{equation}
Thus, for any $x\in\mathbb{R}_+$ we have
$$|T_{g_i}(x)-a_i|\le|T_{g_i}(x)-T_{g_i}(1)|+|T_{g_i}(1)-a_i|
\le\frac1i|x-1|+\frac1i\to 0.$$
as $i\to\infty$. This implies \eqref{approxeq}.
\end{proof}

\begin{remark}\label{rem25}
Observe that if $\frac1n\sum_{i=1}^na_i\to a$ for some infinite
sequence $(a_i)_{i=1}^\infty$ and $\tilde{a}_k$ is a sequence obtained
from $(a_i)_{i=1}^\infty$ by inserting additional elements from a
bounded sequence at sufficiently rarefied times, then
$\frac1n\sum_{i=1}^n\tilde{a}_i\to a$. Therefore the sequence $g_i$ in
the statement of Theorem~\ref{approxim} can be chosen so as to exhaust the
semigroup $G$.
\end{remark}

Note that the proof of Theorem~\ref{approxim} relies strongly on
the contraction property \eqref{contraction} of the action of $G$ on
$\mathbb{R}_+$. The following example shows that without this property
the situation may be different.

\begin{example}
Let $T_r(x)= x+r$ and consider the semigroup $H$ of rational
translations $\{T_r(x),\ r\in\mathbb{Q}_+\}$ acting on $\RR$ (with
$\infty+r=\infty$). For any probability measure $\mu$ concentrated on
$\mathbb{R}_+$ and any infinite sequence of elements $h_i\in H$ there
exists at most one $x\in\mathbb{R}_+$ for which
\begin{equation}\label{example}
\frac1n\sum_{i=1}^n\delta_{T_{h_i}(x)}\stackrel{*}{\longrightarrow}\mu.
\end{equation}
\end{example}

\begin{proof}
Assume there exist $x< y\in\mathbb{R}_+$ such that
$$\lim_{n\to\infty}\frac1n\sum_{i=1}^n\delta_{T_{h_i}(x)}
=\lim_{n\to\infty}\frac1n\sum_{i=1}^n\delta_{T_{h_i}(y)}=\mu.$$
Since $T_{h_i}y=T_{h_i}x+(y-x)$, the sequence \eqref{example} with $x$
replaced by $y$ converges to the measure $\mu$ shifted by $(y-x)$, so
$\mu$ is invariant with respect to the shift by $(y-x)$. However,
there exists an interval $I$ of positive measure $\mu$ and length less
than $(y-x)$. Then $\mu(I+k(y-x))=\mu(I)$ for all $k\ge 0$, so
$\mu(\R_+)=\infty$. This is a contradiction, as $\mu$ is a probability
measure.
\end{proof}


\section{Existence and uniqueness of invariant measures}\label{exuni}


In Ergodic Theory one usually considers invariant measures. As we
explained in the Introduction, requiring the measure to be invariant for
each $T_\omega$ is too strong. Therefore,
we adopt Definition~\ref{definv1}, which takes into account the
particular structure of $G$ as a semigroup with two generators.

Note that in the theory of iterations of random maps (see, e.g.,
\cite{Kif-Liu}) one often encounters a definition of invariant measure
similar to our Definition~\ref{definv1}. This is also the situation in
the theory of Iterated Function Systems (IFS; see, e.g., \cite{BE},
\cite{DF}, \cite{NSB} or \cite{PSS}), except that the
maps $T_{\omega_i}$ are applied in the reverse order. Hence, instead
of considering the measure $\nu$, we should, in the IFS case, consider
the ``dual'' measure. By a dual measure we mean the projection to the
space of measures on the negative coordinates of the natural extension
of $\nu$ as in Section~\ref{secinv}. When $\nu$ is a
Bernoulli measure, then it is naturally ``isomorphic'' to its dual.

Often, in the theory of random iterations and in IFS only Bernoulli
measures $\nu$ are considered. While this approach can be easily
justified, from a purely mathematical point of view all non trivial
ergodic invariant probability measures $\nu$ are equally interesting
(this approach is also adopted for instance in \cite{PSS}).
This motivates our Definition~\ref{definv1}.

We want to show that measures invariant for $(T_0,T_1,\nu)$ are unique
if they exist. We start by proving that our system is contracting in
the sense revealed by the following lemma.

\begin{lemma}\label{lemcon}
Fix $x,y\in\mathbb{R}_+$ and $\bar\omega=(\omega_0,\omega_1,\dots)\in
\Sigma$. Set $x_0=x$, $y_0=y$ and by induction
$x_{n+1}=T_{\omega_n}(x_n)$ and $y_{n+1}=T_{\omega_n}(y_n)$. Then
\begin{equation}\label{lc}
\lim_{n\to\infty}\left(\frac{1}{(1+x_n)}-\frac{1}{(1+y_n)}\right)=0.
\end{equation}
\end{lemma}

\begin{proof}
We consider another system smoothly conjugated to $(T_0,T_1)$. This
conjugacy gives us a new metric on the positive real half-line,
which enables us to control the estimates related to the
contracting properties of the system more efficiently.

The conjugacy is via the logarithmic function: instead of
$x\in(0,\infty)$ we consider $t=\ln x\in\mathbb{R}$. Let $T_0$ and
$T_1$ be conjugate to $L_0$ and $L_1$ respectively. Then $L_0(t)=
\ln(ae^t)=t+\ln a$ and $L_1(t)=\ln(be^t+1)$. We get $L_0'(t)=1$ and
$L_1'(t)=be^t/(be^t+1)<1$, so both maps are Lipschitz continuous with
constant 1. This shows that $|\ln x_{n+1}-\ln y_{n+1}|\le|\ln x_n-\ln
y_n|$ for all $n$. In particular, if one of the sequences $x_n$, $y_n$ converges to
0 or $\infty$, so does the other one.

Assume that$x_n$ and $y_n$ do not converge to zero or infinity. Then there is
$M$ such that for infinitely many $n$'s we have $\ln x_n,\ln y_n\le M$
and $\omega_n=1$ (we use here the assumption that $a<1$). For each
such $n$ we get
$$|\ln x_{n+1}-\ln y_{n+1}|\le|\ln x_n-\ln
y_n|\cdot\frac{be^M}{be^M+1}.$$
Since $be^M/(be^M+1)<1$, this proves that $\lim_{n\to\infty}|\ln
x_n-\ln y_n|=0$. In order to get $1/(1+x)$ from $\ln x$ we apply the
map $t\mapsto 1/(1+e^t)$. The derivative of this map is
$-e^t/(1+e^t)^2$. However, since $(1+e^t)^2-4e^t=(1-e^t)^2\ge 0$, the
absolute value of this derivative is not larger than $1/4$. Therefore
$|1/(1+x_n)-1/(1+y_n)|\le\frac14|\ln x_n-\ln y_n|$, so \eqref{lc}
holds.

If both $x_n$ and $y_n$ converge to 0 then $1/(1+x_n)$ and $1/(1+y_n)$
go to 1, so \eqref{lc} holds.
If both $x_n$ and $y_n$ converge to $\infty$ then $1/(1+x_n)$ and
$1/(1+y_n)$ go to 0, so \eqref{lc} also holds.

Since we started by taking logarithms of $x$ and $y$, there remains
the case when one of these numbers is 0. However, then either all
$\omega_n$'s are 0, and both $x_n$ and $y_n$ go to zero, or there is
$n$ such $\omega_n=1$ and we can apply the main part of the proof to
$x_{n+1}$ and $y_{n+1}$ instead of $x$ and $y$. This completes the
proof.
\end{proof}

Now we can prove the following theorem.

\begin{theorem}\label{unique}
For any ergodic $\sigma$-invariant probability measure $\nu$ there is
at most one $\tau$-invariant probability measure $\kappa$ on
$X=\Sigma\times\RR$, other than
$\nu\times\delta_\infty$, whose projection to $\mathcal{N}$ is $\nu$.
This measure $\kappa$, if it exists, is ergodic.
\end{theorem}

\begin{proof}
If such a measure $\kappa$ exists and is not ergodic, each of its
ergodic components is also $\tau$-invariant and projects to $\nu$.
Therefore it is enough to prove that there is at most one ergodic
measure $\kappa$ satisfying the assertions of the theorem.

Suppose that there exist two such ergodic measures, $\kappa_1$ and $\kappa_2$. For
$i=1,2$ and $\nu$-almost every $\bar\omega\in\Sigma$, there are
generic points of $\kappa_i$ in $\{\bar\omega\}\times\RR$. Take such
$\bar\omega\in\Sigma$ and $x_1,x_2\in\mathbb{R}_+$ such that
$z_i=(\bar\omega,x_i)$ is generic for $\kappa_i$, $i=1,2$ (note that
$x_i\ne\infty$, since $\kappa_i\ne\nu\times\delta_\infty$).

The metric $d(x,y)=|1/(1+x)-1/(1+y)|$ is compatible with the topology
of the compactification $\RR$ of $\mathbb{R}_+$, so we may use it for
measuring of the distance between two points of $X$ having the same
projection to $\Sigma$. By Lemma~\ref{lemcon} we have
$\displaystyle{\lim_{n\to\infty}d(\tau^n(z_1), \tau^n(z_2))}=0$ and
thus $\displaystyle{\lim_{n\to\infty}
|\phi(\tau^n(z_1))-\phi(\tau^n(z_2))|}=0$ for any continuous real
function $\phi$ on $X$. Since $z_i$ is generic for $\kappa_i$
($i=1,2$), this implies that
$\int\phi\,d\kappa_1=\int\phi\,d\kappa_2$. Since this holds for all
$\phi$, we get $\kappa_1=\kappa_2$.
\end{proof}

\begin{corollary}\label{uniquecoro}
For any ergodic $\sigma$-invariant probability measure $\nu$, there is
at most one probability measure invariant for $(T_0,T_1,\nu)$.
\end{corollary}

We now investigate when such a $\tau$-invariant measure exists.

\begin{theorem}\label{existence}
For any ergodic $\sigma$-invariant probability measure $\nu$ there
exists a probability measure invariant for $(T_0,T_1,\nu)$ if and only
if the Lyapunov exponent $\lambda(\nu)$ is negative. If such a measure
$\mu$ exists, then it is unique and the sequence of measures
\begin{equation}\label{kappa}
\kappa_n=\frac1n\sum_{k=0}^{n-1}\delta_{\tau^k(z)},
\end{equation}
where $z=(\bar\omega,x)$, converges for $\nu$-almost every
$\bar\omega\in\Sigma$ and every $x\in\R_+$ to a unique
$\tau$-invariant measure $\kappa$ whose projection to $\mathcal{M}$ is
$\mu$.
\end{theorem}

\begin{proof}
Suppose first that the Lyapunov exponent $\lambda(\nu)$ of $\nu$ is
negative, but there does not exist any probability measure invariant
for $(T_0,T_1,\nu)$. Let $\bar\omega\in\Sigma$ be a generic point of
$\nu$, and set $z=(\bar\omega,1)\in X$. Suppose that the sequence of measures
$(\kappa_n)_{n\in\N}$, given by \eqref{kappa} with $x=1$, does not converge to
$\nu\times\delta_\infty$ (in the weak-* topology). Then there is a
subsequence of the sequence $(\kappa_n)_{n=1}^\infty$ which converges
to some other measure $\kappa$. By the standard argument (as in the
standard proof of the Krylov-Bogolubov Theorem, see, e.g.\
\cite{Walters}, Theorem~6.9), the measure $\kappa$ is $\tau$-invariant. Since
$\bar\omega$ is a generic point of $\nu$, the projection of $\kappa$
to $\mathcal{N}$ is equal to $\nu$. Note that $\Sigma\times\{\infty\}$
is $\tau$-invariant and $\nu\times\delta_\infty$ is the only invariant
measure concentrated on this set and whose projection to $\mathcal{N}$
is $\nu$. Therefore our measure $\kappa$ has to be an affine
combination of $\nu\times\delta_\infty$ and $\kappa'$, where $\kappa'$
is an $\tau$-invariant probability measure concentrated on
$\Sigma\times \R_+$ and its projection to $\mathcal{N}$ is $\nu$. In
this affine combination the coefficient of $\kappa'$ is non-zero since
we assumed that $\kappa\ne\nu\times\delta_\infty$.
Hence the projection of $\kappa'$ to
$\overline{\mathcal{M}}$ is a probability measure invariant for
$(T_0,T_1,\nu)$. This contradicts the assumption that such a measure
does not exist, thus proving that the sequence
$(\kappa_n)_{n=1}^\infty$ must converge to $\nu\times\delta_\infty$.

Since the Lyapunov exponent of $\nu$ is negative, there exists
$\eta>0$ such that $\nu(C_0)\ln a+\nu(C_1)\ln (b+\eta)<0$. Set
$M=\max(1/a,1/\eta)>1$. For every positive integer $n$ define
$$\begin{array}{rl}
A_n=&\{k\in\{0,1,\dots,n-1\}\,:\,\pi_2(\tau^k(z))\ge M\ \mbox{ and }\
\pi_1(\tau^k(z))=0\},\\
B_n=&\{k\in\{0,1,\dots,n-1\}\,:\,\pi_2(\tau^k(z))\ge M\ \mbox{ and }\
\pi_1(\tau^k(z))=1\},\\
D_n=&\{k\in\{0,1,\dots,n-1\}\,:\,\pi_2(\tau^k(z))<M\}.\\
\end{array}$$
These sets form a partition of $\{0,1,\dots,n-1\}$. Now fix $n$ and
set $\alpha_k=\ln\pi_2(\tau^k(z))$ if $k\in A_n\cup B_n$ and
$\alpha_k=0$ if $k\in D_n$. Let us estimate $\alpha_{k+1}$ from above.

If $k\in A_n$ then either $\alpha_{k+1}=\alpha_k+\ln a$ or
$\alpha_{k+1}=0$. In the second case, since $\alpha_k\ge \ln M$ and
$M\ge 1/a$, we get $\alpha_{k+1}\le\alpha_k+\ln a$; clearly this also
holds in the first case.

If $k\in B_n$ then, if $x=\pi_2(\tau^k(z))$, we have $\alpha_k=\ln x$
and $\alpha_{k+1}=\ln (bx+1)=\ln x+\ln (b+1/x)$. Since $x\ge M\ge
1/\eta$, we get $\alpha_{k+1}(x)\le\alpha_k+\ln (b+\eta)$.

If $k\in D_n$ then, with the same notation, $x<M$, so $\alpha_k=0$ and
$\alpha_{k+1}\le\ln (bx+1)<\ln (bM+1)$. Therefore
$\alpha_{k+1}(x)\le\alpha_k+\ln (bM+1)$.

Summing everything and taking into account that $\alpha_0=0$ (because
$M> 1$), we get
$$\sum_{k=0}^{n-1}\alpha_k=\sum_{k=0}^{n-2}\alpha_{k+1}
\le\sum_{k=0}^{n-1}\alpha_k+|A_n|\cdot\ln a+|B_n|\cdot\ln (b+\eta)
+|D_n|\cdot\ln (bM+1).$$
Therefore
$$\frac{|A_n|}{n}\ln a+\frac{|B_n|}{n}\ln (b+\eta)
+\frac{|D_n|}{n}\ln (bM+1)\ge 0.$$
However, as $n\to\infty$, $\kappa_n$ goes to $\nu\times\delta_\infty$,
and hence, $|A_n|/n$ goes to $\nu(C_0)$, $|B_n|/n$ goes to $\nu(C_1)$,
and $|D_n|/n$ goes to 0. Therefore we get
$$\nu(C_0)\ln a+\nu(C_1)\ln (b+\eta)\ge 0,$$
contrary to our choice of $\eta$. This proves that if the Lyapunov
exponent of $\nu$ is negative then there exists a probability measure
$\mu$ invariant for $(T_0,T_1,\nu)$.

Uniqueness of $\mu$ follows from Corollary~\ref{uniquecoro} and
uniqueness of $\kappa$ from Theorem~\ref{unique}. This shows also that
$\kappa$ is ergodic.

If $z$ is a generic point of $\kappa$, the measures $\kappa_n$, given
by \eqref{kappa}, converge to $\kappa$. By Lemma~\ref{lemcon}, if this
holds for some $z=(\bar\omega,x)$, it holds also for
$z'=(\bar\omega,x')$ for every $x'\in\R_+$. Therefore it holds for
$\nu$-almost every $\bar\omega\in\Sigma$ and every $x\in\R_+$.

Suppose now that the Lyapunov exponent of $\nu$ is non-negative, but
there exists a probability measure invariant for $(T_0,T_1,\nu)$. This
measure is the projection to $\overline{\mathcal{M}}$ of some
$\tau$-invariant measure $\kappa$ on $X$ whose projection to
$\mathcal{N}$ is $\nu$. In particular, $\kappa$ is concentrated on
$\Sigma\times \R_+$. Let $z'=(\bar\omega,y)\in X$ be a generic point
of $\kappa$.

Our assumptions imply $\nu(C_1)>0$, since $\nu(C_1)=0$ would imply
$\lambda(\nu)=\ln a<0$ leading to a contradiction. So there exists
$M'>0$ with $\kappa(\Sigma\times\{M'\})=0$ and $\kappa(C_1\times
[0,M'])>0$, which means that we can choose an arbitrarily large $n$
such that $\pi_2(\tau^{n-1}(z'))\le M'$.

Now we make a construction very similar to the one from the proof of
the first implication. The definitions of the sets $A'_n$, $B'_n$ and
$D'_n$ and the functions $\alpha'_k$ will be very similar. Namely, fix
some $M'$ and $n$ as above and set
$$\begin{array}{rl}
A'_n=&\{k\in\{0,1,\dots,n-1\}\,:\,\pi_1(\tau^k(z'))=0\},\\
B'_n=&\{k\in\{0,1,\dots,n-1\}\,:\,\pi_2(\tau^k(z'))\le M'
\ \mbox{ and }\ \pi_1(\tau^k(z'))=1\},\\
D'_n=&\{k\in\{0,1,\dots,n-1\}\,:\,\pi_2(\tau^k(z'))>M'
\ \mbox{ and }\ \pi_1(\tau^k(z'))=1\}.
\end{array}$$
Again, those sets form a partition of $\{0,1,\dots,n-1\}$.

Set $\alpha'_k=\ln\pi_2(\tau^k(z'))$ for all $k\in\{0,1,\dots, n-1\}$.
Let us estimate $\alpha'_{k+1}$ from below.

If $k\in A'_n$ then $\alpha'_{k+1}=\alpha'_k+\ln a$.

If $k\in B'_n$ then, if $x=\pi_2(\tau^k(z))$, we have $\alpha'_k=\ln
x$ and $\alpha'_{k+1}=\ln (bx+1)=\ln x+\ln
(b+1/x)\ge\alpha'_k+\ln(b+1/M')$.

If $k\in D'_n$ then, with the same notation, $\alpha'_{k+1}=\ln
(bx+1)>\ln bx=\alpha'_k+\ln b$.

Summing everything we get
\begin{eqnarray*}
\sum_{k=0}^{n-1} \alpha'_k &=&
\alpha'_0+\sum_{k=0}^{n-2}\alpha'_{k+1}\\
&\ge& \sum_{k=0}^{n-1}\alpha'_k+\alpha'_0-\alpha'_{n-1}
+|A'_n|\cdot\ln a+ |B'_n|\cdot\ln(b+1/M')+|D'_n|\cdot\ln b.\\
\end{eqnarray*}
Since $\alpha'_0=\ln y$ and $\alpha'_{n-1}\le\ln M'$ by our choice of
$n$, we get
$$\frac{|A'_n|}{n}\ln a+\frac{|B'_n|}{n}\ln(b+1/M')
+\frac{|D'_n|}{n}\ln b+\frac{\ln y-\ln M'}{n}\le 0.$$
As $n\to\infty$ (along a suitable subsequence), $|A'_n|/n$ goes to
$\nu(C_0)$, $|B'_n|/n$ goes to $\kappa(C_1\times[0,M'])$, $|D'_n|/n$
goes to $\kappa(C_1\times(M',\infty))$ and $(\ln y-\ln M')/n$ goes to
0, we get
$$\nu(C_0)\ln a+\kappa(C_1\times[0,M'])\ln(b+1/M')
+\kappa(C_1\times(M',\infty))\ln b\le 0.$$
By assumption $\lambda(\nu)=\nu(C_0)\ln a+\nu(C_1)\ln b\ge 0$, hence
$$0\le\nu(C_0)\ln a+\Bigl(\kappa(C_1\times[0,M'])
+\kappa(C_1\times (M',\infty)\Bigr)\ln b$$
and the last two equations yield
$$\kappa(C_1\times[0,M'])\ln(b+1/M')\le\kappa(C_1\times[0,M'])\ln b.$$
This contradicts our choice of $M'$, thus completing the proof.
\end{proof}


\section{Ergodic properties of invariant measures}\label{erg}


Recall that the elements of the semigroup $G$ are 0-1 sequences
$\omega=(\omega_0,\dots,\omega_{n-1})$, $n\in\NN$. Denote the cylinder
consisting of all elements of $\Sigma$ whose first $n$ positions
coincide with $\omega$ by $C_\omega$ and write
$x_n=x(\omega)=T_{\omega_{n-1}}\circ \cdots\circ T_{\omega_0}(x)$.
For $\bar\omega\in\Sigma$ denote the cylinder whose first $n$
positions agree with the first $n$ positions of $\bar\omega$ by
$C^n_{\bar\omega}$ and write
$x_n=x_n(\bar\omega)=T_{\omega_{n-1}}\circ \cdots\circ
T_{\omega_0}(x).$

In the case of a (semi)group with a single generator $T$, Birkhoff's
Ergodic Theorem can be restated in the following way. If a measure
$\mu$ is invariant and ergodic, then for $\mu$ almost every point $x$
we have
$$\frac1n\sum_{i=0}^{n-1}T_*^i(\delta_x)\stackrel{*}
{\longrightarrow}\mu.$$
We would like to obtain similar results for our semigroup $G$.
What was a single element $T^i$ in the semigroup with one generator now
becomes a whole ``sphere'' $\Sigma_i$ and thus, before averaging with
respect to $i$, one should also take the average over this sphere. Luckily, it turns
out that after averaging over the sphere using appropriate weights, we
do not need to average over $i$.

We start with a simple lemma. For any $x,y\in\R_+$ write
$$d(x,y)=\left|\frac{1}{(1+x_n)}-\frac{1}{(1+y_n)}\right|$$
and $\psi_n(\bar\omega)=d(x_n(\bar\omega),y_n(\bar\omega)).$

\begin{lemma}\label{eta}
For any probability measure $\nu$ on $\Sigma$, any $\eta>0$ and any
$x, y\in\R_+$ we have
$$\lim_{n\to\infty}\sum_{\{\bar\omega\,:\,\psi_n(\bar\omega)
\ge\eta\}}\nu(C^n_{\bar\omega})=\lim_{n\to\infty}\nu\bigl(
\{\bar\omega\in\Sigma\,:\,\psi_n(\bar\omega)\ge\eta\}\bigr)=0.$$
\end{lemma}

\begin{proof}
By Lemma~\ref{lemcon}, the sequence of functions $\psi_n$ converges to
0 for any $x, y\in\R_+$, so it converges to 0 in measure.
\end{proof}

\begin{theorem}\label{balls}
Assume that $\nu$ is an ergodic $\sigma$-invariant probability measure
with $\lambda(\nu)<0$. Then for every $x\in\R_+$ the sequence of
measures
$$\mu_{n,x}=\sum_{\omega\in \Sigma_n}\nu(C_\omega)\delta_{T_\omega(x)}$$
converges (in the weak-* topology) to the $(T_0,T_1,\nu)$-invariant
measure $\mu$.
\end{theorem}

\begin{proof}
Given $\eps>0$, we choose $M$ such that $M>x$ and $\mu([0,M])>1-\eps$.
By Lemma~\ref{eta}, for every $\eta>0$ there exists $N$ such that if
$n\ge N$ then
$$\nu\bigl(\{\bar\omega\in\Sigma\,:\,
d\,(0_n(\bar\omega), M_n(\bar\omega))\ge\eta\}\bigr)<\eps.$$
In other words, if
$$B_n=\bigcup\{C_\omega\colon\omega\in \Sigma_n, d(T_\omega
0,T_\omega M)<\eta\}$$
then $\nu(B_n)>1-\eps$ for every $n\ge N$. Therefore,
\begin{equation}\label{e1}
\kappa(B_n\times[0,M])>1-2\eps.
\end{equation}
For a continuous function $\phi:\RR\to\R$, let $\Phi=\phi\circ\pi_2$. Then
$$\int\phi\;d\mu=\int\Phi\;d\kappa.$$
Set $\kappa_n=\tau^n_*\left(\kappa|_{B_n\times[0,M]}\right)$. Since
$\kappa$ is invariant, we have by \eqref{e1}
$$\left|\int\Phi\;d\kappa_n-\int\Phi\;d\kappa\right|
=\left|\int\Phi\;d\left(\tau^n_*\left(\kappa|_{(\Sigma\times
\R_+)\setminus(B_n\times[0,M])}\right)\right)\right|
\le\|\phi\|\cdot 2\eps.$$
If we denote by $\sum{}'$ the sum over those $\omega\in \Sigma_n$ for
which $C_\omega\subset B_n$, and by $\sum{}''$ the sum over the rest
of $\omega\in \Sigma_n$, we get
$$\int\Phi\;d\kappa_n=\int_{B_n\times[0,M]}\Phi\circ \tau^n\;d\kappa
=\sum{}'\ \int_{C_\omega\times[0,M]}\Phi\circ \tau^n\;d\kappa$$
and
$$\int\phi\;d\mu_{n,x}=\sum_{\omega\in\Sigma_n}\nu(C_\omega)\;
\phi(T_\omega x)=\sum{}'\ \nu(C_\omega)\;\phi(T_\omega x)
+\sum{}''\ \nu(C_\omega)\;\phi(T_\omega x),$$
where
$$\left|\sum{}''\ \nu(C_\omega)\;\phi(T_\omega x)\right|
\le\|\phi\|\cdot\sum{}''\ \nu(C_\omega)\le\|\phi\|\cdot\eps.$$
Therefore
\begin{multline}\label{e7}
\left|\int\phi\;d\mu-\int\phi\;d\mu_{n,x}\right|
=\left|\int\Phi\;d\kappa-\int\phi\;d\mu_{n,x}\right|\\
\le 3\ \eps\cdot\|\phi\|+\sum{}'\
\left|\int_{C_\omega\times[0,M]}\Phi\circ\tau^n\;d\kappa
-\nu(C_\omega)\;\phi(T_\omega x)\right|.
\end{multline}
Since $\phi$ is uniformly continuous, for every $\eps>0$ there exists
$\eta>0$ such that, if $d(x, y)<\eta$ then $|\phi(x)-\phi(y)|<\eps$.
By the definition of $B_n$, if $n\ge N$, then for any $\bar\omega'\in
C_\omega\in B_n$ and any $x, y\in [0,M]$ we have
$$|(\Phi\circ \tau^n)(\bar\omega',y)-\phi(T_\omega x)|=|\phi(T_\omega
y)- \phi(T_\omega x)|<\eps.$$
We then get for the sum on the right-hand side of (\ref{e7}) the
following estimate:
$$\sum{}'\ \left|\int_{C_\omega\times[0,M]}\Phi\circ \tau^n\;d\kappa
-\nu(C_\omega)\;\phi(T_\omega x)\right|\le
\eps\cdot\sum{}'\ \nu(C_\omega)\le\eps.$$
Thus, for $n\ge N$ we have
$\left|\int\phi\;d\mu-\int\phi\;d\mu_{n,x}\right|\le \eps (3\|\phi\|+
1)$. As $\eps>0$ can be chosen arbitrarily this proves that
$\lim_{n\to\infty}\int\phi\,d\mu_{n,x}=\int\phi\,d\mu$. As this holds
for every $\phi$, the proof of the theorem is complete.
\end{proof}

One can also mimic Birkhoff's Ergodic Theorem in another way.
Namely, choose one sequence $\bar\omega\in\Sigma$ and average
$\delta_{x_i(\bar\omega)}$ where $i$ ranges from $0$ to $n-1$. Note
that this is a similar procedure as in Theorem~\ref{approxim},
but here we use the structure of $G$ given by its generators.
Therefore we get a different result, which is a direct consequence of
Theorem~\ref{existence}.

\begin{theorem}
Assume that $\nu$ is an ergodic $\sigma$-invariant probability measure
with $\lambda(\nu)<0$. Then for
$\nu$-almost every $\bar\omega\in\Sigma$ and every $x\in\R_+$ the
sequence of measures
$$\frac1n\sum_{i=0}^{n-1}\delta_{x_n(\bar\omega)}$$
converges (in the weak-* topology) to the $(T_0,T_1,\nu)$-invariant
measure $\mu$.
\end{theorem}


\section{Absolutely continuous invariant measures}\label{secinv}


In the theory of Dynamical Systems, measures that are absolutely
continuous with respect to the Lebesgue measure are of special
interest. As we mentioned in the Introduction, we do not know whether
certain particular measures are absolutely continuous. Nevertheless,
we can show that for some $\sigma$-invariant measures $\nu$ there
exist $(T_0,T_1,\nu)$-invariant measures that are absolutely
continuous.

For any $\gamma>1$, consider that map defined by
\begin{equation}\label{Random2}
T(x)=
\begin{cases}
T_0^{-1}(x)=\frac{x}{a} & \mbox{ if }\ \frac{\gamma-1}{b}\le x <
\gamma;\\
T_1^{-1}(x)=\frac{x-1}{b} & \mbox{ if }\ \gamma\le
x\le\frac{\gamma}{a};\\
\end{cases}
\end{equation}
and let $I_0=[\frac{\gamma-1}{b},\gamma)$,
$I_1=[\gamma,\frac{\gamma}{a}]$ and
$I=I_0\cup I_1=[\frac{\gamma-1}{b},\frac{\gamma}{a}]$.
We want to find an absolutely continuous $T$-invariant measure and
investigate its properties. For this, we need the following theorem.

\begin{theorem}\label{LY}
Let $p<\gamma<q\in R$ and let $S$ be a map of $J=[p, q]$ to itself such
that both $S|_{[p, \gamma)}$ and $S|_{(\gamma, q]}$ extend to monotone
$C^2$ maps on the closed intervals $[p, \gamma]$ and $[\gamma, q]$
respectively. Assume that there exists $n$ such that $S^n$ is
piecewise expanding. Then there exists a unique probability measure
$\mu$ on $J$ which is invariant for $S$ and absolutely continuous with
respect to the Lebesgue measure. The measure $\mu$ is ergodic and its
support is the union of finitely many intervals, one of which contains
$\gamma$ in its interior.
\end{theorem}

\begin{proof}
By results of Lasota and Yorke \cite{Lasota-Yorke1} and Li and Yorke
\cite{Li-Yorke} (see also \cite{Boya-Gora}) there exist absolutely
continuous probability measures $\mu_1, \ldots, \mu_k$, $k\ge 1$,
invariant for $S^n$, ergodic, having mutually disjoint supports and such that
every absolutely continuous probability measure, invariant for $S^n$,
is a convex combination of $\mu_1, \ldots, \mu_k$. Moreover, the
support of each $\mu_i$ is the union of finitely many intervals, one
of which contains a discontinuity of $S^n$ or of its derivative in its
interior. Note that such a discontinuity is a preimage of $\gamma$
under $S^{j(i)}$ for some $0\le j(i)<n$ and the measure
$S_*^{j(i)}(\mu_i)$ thus contains $\gamma$ in the interior of its
support.

Since $S^n$ commutes with $S$, any measure $S_*(\mu_i)$ is an
absolutely continuous probability measure, invariant and ergodic for $S^n$,
Therefore $S_*$ maps $\{\mu_1, \ldots, \mu_k\}$ to itself. Since each
$\mu_i$ is fixed by $S_*^n$, this map is a bijection. Only one of
these measures may contain $\gamma$ in the interior of its support,
hence all $S_*^{j(i)}(\mu_i)$ must be equal. This measure belongs to
the $S_*$-orbit of each $\mu_i$ and therefore all $\mu_i$ form one
orbit of $S_*$.

Set $\mu=\frac1k\sum_{i=1}^k\mu_i$. All the properties of $\mu$ from
the statement of the theorem follow directly from the properties of
the $\mu_1, \ldots, \mu_k$ described above.
\end{proof}

We get

\begin{corollary}\label{acip}
For any $0<a<1<b$ and $\gamma>1$ there exists an absolutely continuous
$T$-invariant probability measure with the properties listed in
Theorem~\ref{LY}.
\end{corollary}

\begin{proof}
{}From Lemma~\ref{MR33} it follows that there exists $n$ such that
$T^n$ is piecewise expanding. We then apply Theorem~\ref{LY}.
\end{proof}

Let $\rho:X\to X$ be the skew product map
$\rho(\bar\omega,x)=(\sigma(\bar\omega),T_{\omega_0}^{-1}(x)).$

\begin{lemma}\label{lift}
Given a $T$-invariant probability measure $\mu$, there exists
$B\subset X$ and a $\rho$-invariant probability measure $\kappa$ on
$B$ such that $(\pi_2)_*(\kappa)=\mu$ and $(\rho|_B, \kappa)$ is
isomorphic to $(T, \mu)$.
\end{lemma}

\begin{proof}
For $x\in I$ consider the itinerary map $\phi: I\to \Sigma$ with
$\phi(x)=(\omega_0,\omega_1,\ldots)\in\Sigma$, where $\omega_k=i$ if
$T^k(x)\in I_i$ (for $k\ge 0$ and $i\in\{0, 1\}$). This map is a Borel
map and its graph $B=\{(\phi(x),x)\,:\,x\in I\}$ is a Borel set in
$X$, and the projection $\pi_2|_B$ is a bijection between $B$ and $I$.
Thus a measure $\mu$ can be transported to $B$ by the inverse
bijection and we obtain a measure $\kappa$ on $B\subset X$. By
construction $(\pi_2)_*(\kappa)=\mu$. We have
$$\pi_2(\rho(\phi(x), x))=T_{\omega_0}^{-1}(x)=T(\pi_2(\phi(x), x)),$$ so
$(\pi_2)_*$ is an isomorphism between $(B,\rho,\kappa)$ and
$(I,T,\mu)$. In particular, $\kappa$ is $\rho$-invariant.
\end{proof}

Denote by $\Sigma_{ts}$ the space of all infinite two-sided 0-1
sequences
$$\bar\xi=(\dots,\xi_{-2},\xi_{-1},\xi_{0},\xi_{1},\xi_{2},\dots).$$

\begin{lemma}\label{muinvariant}
Given an ergodic $T$-invariant probability measure $\mu$, there exists
an ergodic $\sigma$-invariant probability measure $\nu$ on $\Sigma$
such that $\mu$ is invariant for $(T_0,T_1,\nu)$.
\end{lemma}

\begin{proof}
Let us look at the natural extension (inverse limit)
$(\widehat{X},\widehat\rho,\widehat\kappa)$ of the system
$(X,\rho,\kappa)$ constructed in Lemma~\ref{lift}. The space
$\widehat{X}$ consists of all sequences of pairs
$(\bar\omega^{(i)},y^{(i)})_{i=-\infty}^0$, for which
\begin{equation}\label{natex}
\rho(\bar\omega^{(i)},y^{(i)})=(\bar\omega^{(i+1)},y^{(i+1)})
\end{equation}
($i=-1,-2,\dots$). The map $\widehat\rho$ is given by
\begin{multline}
\widehat\rho(\dots,(\bar\omega^{(-2)},y^{(-2)}),
(\bar\omega^{(-1)},y^{(-1)}), (\bar\omega^{(0)},y^{(0)}))\\
=(\dots,(\bar\omega^{(-1)},y^{(-1)}),(\bar\omega^{(0)},y^{(0)}),
(\sigma(\bar\omega^{(0)}),T_{\omega^{(0)}_0}^{-1}(y^{(0)}))).
\end{multline}

Condition \eqref{natex} is equivalent to the pair of conditions
\begin{equation}\label{natex1}
\sigma(\bar\omega^{(i)})=\bar\omega^{(i+1)}
\end{equation}
and
\begin{equation}\label{natex2}
T_{\omega^{(i)}_0}^{-1}(y^{(i)})=y^{(i+1)}.
\end{equation}
Sequences $(\bar\omega^{(i)})_{i=-\infty}^0$ of elements of $\Sigma$
satisfying \eqref{natex1} can be identified with doubly infinite 0-1
sequences $\bar\xi$ by $\bar\omega^{(i)}=(\xi_i,\xi_{i+1},\dots)$.
Then \eqref{natex2} becomes
\begin{equation*}
T_{\xi_i}^{-1}(y^{(i)})=y^{(i+1)}.
\end{equation*}
This is equivalent to $T_{\xi_i}(y^{(i+1)})=y^{(i)}$. Therefore, once
$\bar\xi$ and $y^{(0)}$ are given, all $y^{(i)}$ are uniquely
determined and become redundant. This means that the space
$\widehat{X}$ can be written as $\Sigma_{ts}\times\RR$, and then
$\widehat\rho$ becomes the skew product
$$\widehat\rho(\bar\xi,y)=(\sigma(\bar\xi),T_{\xi_0}^{-1}(y)).$$
Note that although $T_{\xi_0}^{-1}(y)$ is not defined everywhere, it
is defined $\widehat\kappa$-almost everywhere, and this is sufficient for
our purposes.

The inverse of $\widehat\rho$ is given by
$$\widehat\rho^{\
-1}(\bar\xi,y)=(\sigma^{-1}(\bar\xi),T_{\xi_{-1}}(y)).$$
Therefore $\tau$ is a factor of $\widehat\rho^{\ -1}$ under the
projection $\pi:\widehat{X}\to X$, given by
$$\pi(\bar\xi,y)=((\xi_{-1},\xi_{-2},\dots),y).$$
Hence, $\widetilde\kappa=\pi_*(\widehat\kappa)$ is a $\tau$-invariant
probability measure.

By Lemma~\ref{lift}, $(\pi_2)_*(\kappa)=\mu$. Therefore for every
measurable set $A\subset\R_+$ we have $\mu(A)=\kappa(\Sigma\times A)$.
Clearly, $\kappa(\Sigma\times A)= \widehat\kappa(\Sigma_{ts}\times
A)$. This proves that the projection of $\widehat\kappa$ to the space
of probability measures on the second coordinate of
$\widehat{X}=\Sigma_{ts}\times\RR$ is $\mu$. Hence,
$(\pi_2)_*(\widetilde\kappa)=\mu$. Therefore, if
$\nu=(\pi_1)_*(\widetilde\kappa)$, then $\mu$ is a
$(T_0,T_1,\nu)$-invariant measure.

Since $\mu$ is ergodic for $T$, by Lemma~\ref{lift} $\kappa$ is
ergodic for $\rho$. Therefore $\widehat\kappa$ is ergodic for
$\widehat\rho$, so $\widetilde\kappa$ is ergodic for $\tau$, so $\nu$
is ergodic for $\sigma$.
\end{proof}

We can now state the main result of this section.

\begin{theorem}\label{abscon}
For every $\gamma>1$ there exists a $\sigma$-invariant ergodic
probability measure $\nu_\gamma\in\mathcal{N}$ with negative Lyapunov
exponent, for which the corresponding $(T_0,T_1,\nu_\gamma)$-invariant
measure $\mu_\gamma$ is absolutely continuous with respect to the
Lebesgue measure and the support of $\mu_\gamma$ is the union of
finitely many intervals whose convex hull equals
$[\frac{\gamma-1}{b},\frac{\gamma}{a}]$. In particular,
$\mu_\gamma\neq \mu_\delta$ for $\gamma\neq\delta$, $\gamma,\delta>1$.
\end{theorem}

\begin{proof}
The existence of $\nu_\gamma$ and $\mu_\gamma$ follows from
Corollary~\ref{acip}, Lemma~\ref{muinvariant} and
Theorem~\ref{existence}. Observe that by Theorem~\ref{LY} $\gamma$ is
contained in the interior of the support of $\mu_\gamma$. Since
$T(\gamma)=\frac{\gamma-1}{b}$ and $\lim_{x\to \gamma_-}
T(x)=\frac{\gamma}{a}$, the endpoints of $I$ are contained in the
support of $\mu_\gamma$. As the support of $\mu_\gamma$ is contained
in $I$, its convex hull equals $I$. All other properties of
$\mu_\gamma$ follow from Corollary~\ref{acip}.
\end{proof}


\section{Invariant measures over Bernoulli systems}\label{bern}


Let us see how Definition \ref{definv1} works in the
classical case when the system $(\sigma,\Sigma,\nu)$ is Bernoulli.

\begin{proposition}\label{product}
Assume that $\nu$ is Bernoulli. If a $\tau$-invariant measure $\kappa$
projects to $\nu$, then $\kappa=\nu\times\mu$ for some measure $\mu$.
\end{proposition}

\begin{proof}
We have $\kappa=\kappa_1+\kappa_2$ where $\kappa_1$ is concentrated on
$\Sigma\times\mathbb{R}_+$ and $\kappa_2$ is concentrated on
$\Sigma\times\{\infty\}$ (one of these measures may be zero). Since
$\Sigma\times\{\infty\}$ is $\tau$-invariant, the measure $\kappa_2$
is $\tau$-invariant. Together with the ergodicity of $\nu$ this
implies that $\kappa_2=c\cdot\nu\times\delta_\infty$ for some constant
$c$. Therefore it is enough to consider $\kappa=\kappa_1$. Thus we
assume that $\kappa$ is concentrated on $\Sigma\times\R_+$.

Observe first that the measure $\tau_*(\nu\times\delta_x)$ is a
weighted average of the measures $\nu\times\delta_{T_0x}$ and
$\nu\times\delta_{T_1x}$ (here we use the assumption that $\nu$ is
Bernoulli). Therefore, for any measure $\mu_1$ on $\RR$, the measure
$\tau_*(\nu\times\mu_1)$ is of the form $\nu\times\mu_2$ for some
measure $\mu_2$ on $\RR$.

By Theorem~\ref{existence}, for almost every $\bar\omega\in\Sigma$,
the sequence of measures $\kappa_n$ given by \eqref{kappa} with $x=1$
converges to $\kappa$. Integrating those measures with respect to
$\nu$ and applying the Lebesgue Dominated Convergence Theorem (we can
do it since the weak\nobreakdash-* convergence means convergence of integrals of
any continuous function $\phi:\RR\to\R$ and those integrals are
bounded by the supremum of $|\phi|$) we get that
$$\frac1n\sum_{i=0}^{n-1}\tau_*^i(\nu\times\delta_1)\stackrel{*}
{\longrightarrow}\kappa.$$
By the preceding paragraph, all measures $\tau_*^i(\nu\times\delta_1)$
are of the form $\nu\times\mu_i$, so $\kappa$ is also of the form
$\nu\times\mu$.
\end{proof}

If the system $(\sigma,\Sigma,\nu)$ is $(p,1-p)$ Bernoulli, then the
action of $G$ has a simple interpretation. Namely, we apply $T_0$ or
$T_1$ randomly each time, choosing $T_0$ with probability $p$ and
$T_1$ with probability $1-p$.

In this case we can also look at the invariant measures for
$(T_0,T_1,\nu)$ from another point of view. The maps $T_0,T_1$ induce
operators $(T_0)_*,(T_1)_*$ from $\mathcal{M}$ to itself. Define
$\mathcal{T}_p:\mathcal{M}\to\mathcal{M}$ by $\mathcal{T}_p=p(T_0)_*+(1-p)(T_1)_*$.

\begin{proposition}\label{invber}
Let $\nu$ be the $(p,1-p)$ Bernoulli measure on $\Sigma$ and let
$\mu\in\mathcal{M}$. Then the following are equivalent
\begin{enumerate}
\item the measure $\nu\times\mu$ is $\tau$-invariant;
\item for any measurable set $B\subset\RR$, we have
\begin{equation}\label{mube}
\mu(B)=p\mu(T_0^{-1}(B))+(1-p)\mu(T_1^{-1}(B));
\end{equation}
\item $\mu$ is a fixed point of $\mathcal{T}_p$.
\end{enumerate}
\end{proposition}

\begin{proof}
Denote $C_i=\{\bar\omega\in\Sigma\,:\,\omega_0=i\}$ for $i=0,1$. We
have
$$\tau^{-1}(A\times B)=\big[(\sigma|_{C_0})^{-1}(A)\times T_0^{-1}(B)
\big]\cup \big[(\sigma|_{C_1})^{-1}(A)\times T_1^{-1}(B)\big].$$
Since $\nu\bigl((\sigma|_{C_0})^{-1}(A)\bigr)=p\nu(A)$ and
$\nu\bigl((\sigma|_{C_1})^{-1}(A)\bigr)=(1-p)\nu(A)$, we have
\begin{equation}\label{nucrossmu}
(\nu\times\mu)(\tau^{-1}(A\times B))=
p\nu(A)\mu(T_0^{-1}(B))+(1-p)\nu(A)\mu(T_1^{-1}(B)).
\end{equation}
If  $\nu\times\mu$ is $\tau$-invariant, then by \eqref{nucrossmu}
$$\mu(B)=(\nu\times\mu)(\Sigma\times B)=(\nu\times\mu)\tau^{-1}
(\Sigma\times B)=p\mu(T_0^{-1}(B))+(1-p)\mu(T_1^{-1}(B)),$$
so (2) holds.

On the other hand, multiplying both sides of \eqref{mube} by
$\nu(A)$ and using \eqref{nucrossmu} one sees that (2) implies (1).

Conditions (2) and (3) are equivalent by the definition of $\mathcal{T}_p$.
\end{proof}

Note that by the above proposition and Theorem~5.1 of \cite{DF}
(cf.~\cite{NSB}), our $(T_0, T_1, \nu)$-invariant measure is the same
as the invariant (stationary) measure for the corresponding IFS.

As an immediate consequence of Propositions \ref{product} and
\ref{invber} we get the following theorem.

\begin{theorem}
Assume that $\nu$ is $(p, 1-p)$ Bernoulli. A probability measure $\mu$
supported by $\R_+$ is $(T_0, T_1, \nu)$-invariant if and only if it
is a fixed point of $\mathcal{T}_p$.
\end{theorem}

Another immediate consequence is the following lemma, where we use
property (2) of Proposition~\ref{invber}.

\begin{lemma}\label{measuremult}
Assume that $\nu$ is $(p, 1-p)$ Bernoulli and $\mu$ is the $(T_0, T_1,
\nu)$-invariant measure. Then for every interval $I$ and $i=0,1$ we
have $\mu(T_i(I))\ge q\mu(I)$, where $q=\min(p,1-p)$.
\end{lemma}

{}From this we get, in turn, the next lemma (cf.~\cite{BE},
Theorem~3).

\begin{lemma}\label{support}
Assume that $\nu$ is $(p, 1-p)$ Bernoulli and $\mu$ is the $(T_0, T_1,
\nu)$-invariant measure. Then the measure $\mu$ of every interval $I$
is positive. In other words, the support of $\mu$ is the whole $\RR$.
\end{lemma}

\begin{proof}
There exists an interval $J\subset\R_+$ with $\mu(J)>0$. Let $x$ be
the center of $J$. Let $I\subset\R_+$ be an interval with center $y$.
By Lemmas~\ref{MR33} and~\ref{MR34}, we can find $\omega\in G$ such
that that $T_\omega(x)$ is so close to $y$ and $T_\omega'$ is so small
that $T_\omega(J)\subset I$. Then by Lemma~\ref{measuremult}
$\mu(I)>0$.
\end{proof}

As we mentioned in the Introduction, we do not know whether the assumption that $\nu$
is $(p, 1-p)$ Bernoulli and $\mu$ is a $(T_0, T_1, \nu)$-invariant
measure implies that $\mu$ is absolutely continuous with respect to the
Lebesgue measure. The strongest property of the measure $\mu$ that we
can prove is the following.

\begin{theorem}\label{holderimage}
Assume that $\nu$ is $(p,1-p)$-Bernoulli with $\lambda(\nu)<0$. Then
the $(T_0,T_1,\nu)$-invariant measure $\mu$ on $\R_+$ is the image of
the Lebesgue measure under an increasing map from $[0,1)$ onto $\R_+$,
which is H\"older continuous with the same exponent on each compact
interval.
\end{theorem}

In order to prove this theorem we need the following lemma. For
an interval $I$ we denote its length by $|I|$.

\begin{lemma}\label{technical}
Under the assumptions of Theorem~\ref{holderimage}, if $I\subset\R_+$
is an interval centered at $t$, then
\begin{equation}\label{8}
\log\mu(I)>c_4(t)+c_5\log|I|
\end{equation}
holds with $c_4(t)$ depending continuously on $t$ (but independent of
$I$) and
\begin{equation}\label{9}
c_5=\frac{-\log q\cdot\log\frac{a+b-1}{a}}{\log b\cdot\log\frac{b}
{a+b-1}}>0,
\end{equation}
where $q=\min(p,1-p)$.
\end{lemma}

\begin{proof}
For an interval $I\subset \R_+$, we need to estimate $\mu(I)$ from
below. If we find a sequence $\omega$ of length $n$ such that
$T_\omega([0,1/a])\subset I$, then by Lemma~\ref{measuremult}, $\mu(I)\ge
q^n\mu([0,1/a])$. Note that by Lemma~\ref{support}, $\mu([0,1/a])>0$.
Thus, we will look for preimages of $I$ along a specific branch.

We start by choosing the smallest $k$ such that for the midpoint $t$
of $I$ we have $T_1^{-k}(t)\le 1/a$. Since $T_1(0)=1<1/a$, we have
then $T_1^{-k}(t)\in [0,1/a]$. Moreover, if $t>1/a$ then
$1/a<T_1^{-k+1}(t)<t/b^{k-1}$, so
\begin{equation}\label{k}
k<\frac{\log^+ ta}{\log b}+1,
\end{equation}
where we use the notation $\log^+x=\max(0,\log x)$.

Consider the map $\phi:[0,1/a)\to [0,1/a)$ given by
\begin{equation}\label{cases}
\phi(x)=\begin{cases}
\frac{x}{a}& \text{ if } x<1,\\
\frac{x-1}{b}& \text{ if } x\ge 1.
\end{cases}\end{equation}
(cf.\ formula~\eqref{Random2}).
Starting at the point $T_1^{-k}(t)$ and applying $\phi$ amounts to
choosing further preimages of $t$ along some branch. Our aim is to
find $m$ such that $(\phi^m)'(T_1^{-k}(t))$ is larger than $1/a$
divided by the length of the half of the interval $T_1^{-k}(I)$. Then
the preimage of $I$ under the corresponding branch (of length $k+m$)
will contain $[0,1/a)$.

Since we do not know the location of $T_1^{-k}(t)$, we have to
estimate from below $(\phi^n)'(x)$ for all $x\in[0,1/a)$. We can treat
$\phi$ as a map of the circle with one discontinuity at which the
left-hand limit is larger than the right-hand limit (so its lifting is
an \emph{old heavy} map,\footnote{A map $\Phi:\R\to\R$ is \emph{old} if
$\Phi(x+1)=\Phi(x)+1$ for all $x\in\R$; it is \emph{heavy} if
$\lim_{y\to x^-}\Phi(y)\ge\Phi(x)\ge\lim_{y\to x^+}\Phi(y)$ for all
$x\in\R$.} see \cite{M}). We will compare it with the map
$\psi:[0,1/a)\to [0,1/a)$ given by
$$\psi(x)=\begin{cases}
\frac{1-a}{ab}+\frac{a+b-1}{ab}\;x& \text{ if } x<1,\\
\frac{x-1}{b}& \text{ if } x\ge 1.
\end{cases}$$
This map can be treated as an orientation preserving homeomorphism of
the circle. Let us find its rotation number. We have
$$\frac{1-a}{ab}+\frac{a+b-1}{ab}\left(y-\frac{1}{b-1}\right)=
\frac{a+b-1}{ab}y-\frac{1}{b-1}$$
and
$$\frac{\left(y-\frac{1}{b-1}\right)-1}{b}=\frac{y}{b}-\frac{1}{b-1}\;.$$
This means that the map
$$\xi:\left[0,\frac{1}{a}\right)\to\left[\log\frac{1}{b-1},\log
\left(\frac{1}{a}+\frac{1}{b-1}\right)\right),$$
given by the formula
$$\xi(x)=\log\left(x+\frac{1}{b-1}\right),$$
conjugates $\psi$ with the map
$$\zeta:\left[\log\frac{1}{b-1},\log\left(\frac{1}{a}+\frac{1}{b-1}\right)
\right)\to\left[\log\frac{1}{b-1},\log\left(\frac{1}{a}+\frac{1}{b-1}
\right)\right),$$
given by the formula
$$\zeta(x)=\begin{cases}
x+\log\frac{a+b-1}{ab}&\text{if }\log\frac{1}{b-1}\le
x<\log\frac{b}{b-1},\\
x+\log\frac{1}{b}&\text{if }\log\frac{b}{b-1}\le x<\log\left(\frac{1}{a}+
\frac{1}{b-1}\right).
\end{cases}$$

The map $\zeta$, treated as a circle homeomorphism, is just a
rotation. Its rotation number is
$$\rho(\zeta)=\frac{\log\frac{a+b-1}{ab}}{\log\frac{a+b-1}{a}}\;.$$
The lifting of $\zeta$, with the length of the circle normalized to 1,
is the translation by $\rho(\zeta)$. The $n$-th iterate of this
lifting sends a point $y\in[0,1)$ to $y+n\rho(\zeta)\in[p,p+1)$, where
$p$ equals $\lfloor n\rho(\zeta)\rfloor$ or $\lfloor n\rho(\zeta)+
1\rfloor$. Thus, the same is true for $\psi$ replacing $\zeta$.

Let $\Psi$ be the lifting of $\psi$ and $\Phi$ the lifting of $\phi$,
agreeing with $\Psi$ on the right lap. Then
$\Psi\ge\Phi$. Moreover, $\Psi$ is increasing, and thus
$\Psi^n\ge\Phi^n$ for all $n$ (see \cite{ALM}). Thus if $y\in[0,1)$
then $\Phi^n(y)\in[p',p'+1)$, where $p'\le\lfloor
n\rho(\zeta)+1\rfloor$. This means that if we iterate $\phi$ $n$
times, its right lap is used at most $p'$ times. Therefore the left
lap is used at least $n-p'\ge n-2-n\rho(\zeta)$ times. From
\eqref{cases} we read off the derivatives of both laps of $\phi$, and
we get that
$$(\phi^n)'(x)\ge\left(\frac{1}{a}\right)^{n-2-n\rho(\zeta)}
\cdot\left(\frac{1}{b}\right)^{2+n\rho(\zeta)}=
\left(\frac{a}{b}\right)^2\left[\left(\frac{1}{a}\right)\cdot
\left(\frac{a}{b}\right)^{\rho(\zeta)}\right]^n.$$
Substituting the value of $\rho(\zeta)$, we get after a short
algebraic manipulation on the formulas
\begin{equation}\label{phi}
\log(\phi^n)'(x)\ge 2\log\frac{a}{b}+\frac{\log
b\cdot\log\frac{b}{a+b-1}}{\log\frac{a+b-1}{a}}\;n.
\end{equation}

We want to find $m$ such that
$$(\phi^m)'(x)\ge\frac{1}{a}\cdot\frac{2}{|T_1^{-k}(I)|}=
\frac{2b^k}{a|I|}\;.$$
According to \eqref{k} and \eqref{phi}, this inequality will be
satisfied if
$$2\log\frac{a}{b}+\frac{\log b\cdot\log\frac{b}{a+b-1}}{\log
\frac{a+b-1}{a}}\;m\ge\log\frac{2b^{(\log^+ta)/\log b+1}}{a|I|}\;,$$
that is,
$$m\ge\frac{\log\frac{2b^{(\log^+ta)/\log b+3}}{a^3|I|}\cdot\log
\frac{a+b-1}{a}}{\log b\cdot\log\frac{b}{a+b-1}}\;.$$
Therefore we can find required integer $m$ satisfying
$$m<\frac{\log\frac{2b^{(\log^+ta)/\log b+3}}{a^3|I|}\cdot\log
\frac{a+b-1}{a}}{\log b\cdot\log\frac{b}{a+b-1}}+1.$$
We can rewrite this inequality as
\begin{equation}\label{m}
m<c_1(t)-c_2\log|I|,
\end{equation}
where
\begin{equation}\label{5}
c_1(t)=\frac{\log\frac{2b^{(\log^+ta)/\log b+3}}{a^3}\cdot\log
\frac{a+b-1}{a}}{\log b\cdot\log\frac{b}{a+b-1}}+1
\end{equation}
and
\begin{equation}\label{6}
c_2=\frac{\log\frac{a+b-1}{a}}{\log b\cdot\log\frac{b}{a+b-1}}.
\end{equation}

Note that all logarithms appearing in the formulas for $c_1(t)$ and
$c_2$ are positive, so $c_1(t)$ and $c_2$ are positive.

In such a way we found a 0-1 sequence $\omega$ of length $k+m$, such
that $T_\omega([0,1/a])\subset I$, and hence $\mu(I)\ge
q^{(k+m)}\mu([0,1/a])$. Moreover, by \eqref{k} and \eqref{m}, we have
the following estimate for $k+m$:
\begin{equation}\label{7}
k+m<\frac{\log^+ta}{\log b}+1+c_1(t)-c_2\log|I|.
\end{equation}
We can rewrite it as
$$k+m<c_3(t)-c_2\log|I|.$$
According to formulas \eqref{m}-\eqref{7}, the dependence of $c_3(t)$
on $t$ is continuous, and $c_2$ is independent of $t$. We get
$$\log\mu(I)>(c_3(t)-c_2\log|I|)\log q+\log\mu([0,1/a]).$$
Therefore \eqref{8} holds with $c_5$ given by \eqref{9} and $c_4(t)$ depending
continuously on $t$.
\end{proof}

\begin{proof}[Proof of Theorem~\ref{holderimage}]
Consider the map $H:[0,1)\to\R_+$, such that $H(x)=y$ if and only if
$\mu([0,y])=x$. By Lemma~\ref{technical}, and since $\mu(\R_+)=1$,
this map is well defined. By the definition, the measure $\mu$ is the
image of the Lebesgue measure under $H$. In particular, for an
interval $J\subset[0,1)$ we have $|J|=\mu(H(J))$, so \eqref{8} gives
us
$$|J|>\exp(c_4(t))\cdot|H(J)|^{c_5},$$
where $t$ is a certain point of $H(J)$. This shows that $H$ is
H\"older continuous with exponent $1/c_5$ on each compact
interval.
\end{proof}


\end{document}